\documentclass[a4paper,12pt]{article}
\usepackage{graphicx}
\usepackage{amsfonts}
\usepackage{epsfig}
\usepackage{subfigure}
\usepackage{amsmath}

\newtheorem{thm}{Theorem}[section]

\newtheorem{lem}[thm]{Lemma}
\newtheorem{rem}[thm]{Remark}
\newtheorem{cor}[thm]{Corollary}

\newtheorem{exm}[thm]{Example}
\newcommand{\thmref}[1]{Theorem~{\rm \ref{#1}}}
\newcommand{\lemref}[1]{Lemma~{\rm \ref{#1}}}
\newcommand{\corref}[1]{Corollary~{\rm \ref{#1}}}

\newcommand{\remref}[1]{Remark~{\rm \ref{#1}}}
\newcommand{\exmref}[1]{Example~{\rm \ref{#1}}}

\newcommand{\beq}[1]{\begin{equation} \refstepcounter{neweqn} \label{#1}}
\newcommand{\eeq}{\end{equation}}
\newcommand{\bed}{\begin{displaymath}}
\newcommand{\eed}{\end{displaymath}}
\newcommand{\ben}{\begin{eqnarray*}}
\newcommand{\een}{\end{eqnarray*}}

\newcommand{\bedd}{\bed\begin{array}{l}}
\newcommand{\eedd}{\end{array}\eed}
\newcommand{\al}{\alpha}

\newcommand{\nd}{\noindent}

\newcommand{\wdt}{\widetilde}

\newcommand{\cd}{(\cdot)}

\def\({\left(}
\def\){\right)}
\def\one{{\hbox{1{\kern -0.35em}1}}}
\newcommand{\bdd}{\hspace*{-0.08in}{\bf.}\hspace*{0.05in}}
\newcommand{\bea}{\bed\begin{array}{rl}}
\newcommand{\eea}{\end{array}\eed}
\newcommand{\ad}{&\!\!\!\disp}

\newcommand{\barray}{\begin{array}{ll}}
\newcommand{\earray}{\end{array}}

\def\para#1{\vskip 0.4\baselineskip\noindent{\bf #1}}
\def\qed{$\qquad \Box$}

\newcommand{\e}{\varepsilon}

\newcommand{\disp}{\displaystyle}
\newcommand{\T}{{\cal T}}

\title{An
$\e$-uniform Finite Element Method for Singularly Perturbed
Boundary Value Problems}
\author{ Q. S. Song\thanks{Department of Mathematics, Wayne State
University, Detroit, MI 48202, {\tt song@math.wayne.edu}. Research
of this authors was supported in part by a WSU graduate research
assistantship.} \and G. Yin\thanks{Department of Mathematics, Wayne
State University, Detroit, MI 48202, {\tt gyin@math.wayne.edu}.
Research of this author was supported in part by the National
Science Foundation, and in part by Wayne State University Research
Enhancement Program.} \and Z. Zhang\thanks{Department of
Mathematics, Wayne State University, Detroit, MI 48202, {\tt
zzhang@math.wayne.edu}. Research of this author was supported in
part by the National Science Foundation, and in part by Michigan
Life Science Corridor.} }

\begin{document}
\maketitle

\begin{abstract}
This work develops an $\e$-uniform finite element method for
singularly perturbed boundary value problems. A surprising and remarkable observation is illustrated: By moving one node arbitrarily in between its adjacent nodes, the new finite element solution always intersect with original one at fixed point. Using this fact, an effective $\e$-uniform approximation out of boundary is proposed by adding one point only in the grid that contains boundary layer. The thickness of boundary layer is not necessary to be
known from {\it priori estimation}. Numerical results
are carried out and  compared to Shishkin mesh for
demonstration purpose.

\vskip 0.2 in \nd{\bf Key Words.} finite element method, singular
perturbation, $\e$-uniform approximation,
layer-adapted mesh, Shishkin mesh.

\vskip 0.2 in \nd{\bf Mathematics Subject Classification.}

\vskip 0.2 in \nd{\bf Brief Title.} An $\e$-uniform Approximation of
Singularly Perturbed BVP

\end{abstract}

\newpage
\setlength{\baselineskip}{0.28in}

\section{Introduction}
This paper is concerned with linear Galerkin finite element method
for singularly perturbed boundary value problems (BVPs).
Consider an
one-dimensional BVP problem
\begin{equation} \label{bvp1}
-\e u''- bu' + cu = f, \quad u(0) = u(1) = 0, \quad x\in [0,1].
\end{equation}
For simplicity, let $b\le 0, c\ge 0$, and $0<\e \ll 1$ are constant
such that not both $b$ and $c$ are 0.
%with $bc\neq 0$.
If $b>0$, by using substitution $w(x) = u(1-x)$,
%the model meets the assumption of
it reduces to the case with $b\le 0$.
If $b=0, c>0$, equation
(\ref{bvp1}) is said to be a reaction diffusion equation. If
$b<0, c=0$, equation (\ref{bvp1}) is
the so-called convection diffusion
equation. All the results
presented in this paper can be
readily generalized to smooth and
non-vanishing functions of $b(x)$ and $c(x)$.

If the exact solution $u\cd$ of (\ref{bvp1}) is ``bad'' in the sense
that $\|u''\|_\infty$ is not  bounded uniformly in $\e$,  the
standard finite element method (FEM) generates huge errors through
the whole domain. Typically, it is caused by
a small interval of
width $O(\e)$ (called boundary layer), in which $u''$ rapidly changes.

To overcome the difficulties in the singular perturbation, it is
desirable to put more grid points near the boundary layer or stablize the appoximation methods. Streamline diffusion finite element methods (SDFEM), upwinding FEM, Bakhalov grid, Shishkin grid, and many other such schemes are extensively studied in the context of singularly purturbed problems since 1970s, see \cite{Roos1, OR1, Doolan, Miller, Shishkin, Bakhalov}. Among them, Shishkin grid became popular due to its simple structure and high accuracy. The Shishkin mesh was first introduced in finite difference methods and has been discussed in \cite{Miller}; the reader is  referred to a survey article \cite{Roos3} for further details. A typical Shishkin mesh is to construct $n+n$ grid, which
is indeed $n$ uniform grids in boundary layer plus $n$ uniform grids out of boundary layer. By this method, the approximation provides $\e$-uniform accuracy. But they require {\it a priori} estimation in order to determine the thickness of the bounded layer. On the other hand, it makes error analysis more complicated, since the errors from boundary layer affect the solution in the entire domain. Therefore, if an approximation can be stabilized and $\e$-uniform by simply adding one point to original $n$ grid, it deserves to be worked out.

In this work, we focus on FEM solutions of (\ref{bvp1}) by
starting with an interesting observation. Given  a grid $\T^n =
\{0=x_0< \cdots <x_n<x_{n+1}=1\}$, we add $m$ points arbitrarily in
$[x_n, 1]$, denoted by $\T^{n+m}$. Then FEM solutions on $\T^n$ and
$\T^{n+m}$ intersect each other at a fixed point in each interval
of out of boundary layer, that is, the locations of intersection in each interval is indpendent of $m$ and the distribution of added points, see Figure \ref{figcde}-a, Figure \ref{figrde1}-a, and Figure \ref{figgre1}. This directly implies that the accuracy on
those intersections are as good as FEM solutions on the grid $\T^{n+m}$  by choosing $m\to \infty$, denoted by $\T^{n+\infty}$. Provided that the boundary layer is covered by $[x_n,1]$, the above observation gives the start point of $\e$-uniform approximation. In lieu of interpolating these intersections, we present a better way to obtain an
$\e$-uniform approximation. By adding one point $\hat s_1\in
(x_n,1)$ with $|\hat s_1-x_n| = O(\e)$ or  $|\hat s_1-x_n| =
O(\sqrt{\e})$, the interval $[x_n,\hat s_1]$ block the error impact
from boundary layer completely. The theoretical result shows that
the FEM solutions with grid $\T^n\cup\{\hat s_1\}$ in $[0,x_n]$ is
the same as the FEM solution of
\begin{equation} \label{bvp1-1}
-\e w''- bw' + cw = f, \quad w(0) = 0, \  w(\hat s_1) = u(\hat s_1),
\quad x\in [0,\hat s_1],
\end{equation}
where $u(\hat s_1)$ is exact solution at the point of $\hat{s_1}$, and $w''\cd$ is uniformly bounded. This enables us to use all kinds of standard FEM error analysis in $[0,\hat{s_1}]$, no matter how huge errors are generated in $[\hat{s_1},1]$. Therefore,
the FEM errors in $[0,x_n]$ has the accuracy of $\T^{n+1+\infty}$,
which is clearly better than Shishkin mesh $\T^{n+n}$ in both accuracy and computing
cost. Another advantage is: One
need not know the thickness of boundary layer, since $\hat s_1$ is
not necessarily in boundary layer.

The rest of this paper is arranged as follows. Section 2 begins with the
model and notation. Section 3 proceeds with the observation on
intersections of a family of FEM solutions. Section 4 presents
an $\e$-uniform FEM, which can isolate the boundary layer. Some
auxiliary results are included in Section 5. Section 6 displays some
numerical experiment results,
%:
including solutions of convection-diffusion equation,
reaction-diffusion equation,
and Green function.
Finally, we close this
paper with further remarks.

\section{Formulation}

Let $H^1 = \{ v, v' \in L^2\}$, and $H_0^1 = \{ v| v\in H^1 , v(0) =
v(1) = 0\}$. The weak solution of (\ref{bvp1}) is a function $u\in
H_0^1$, satisfying
\begin{equation} \label{wksol1}
a(u,v) = (f,v), \forall v\in H_0^1,
\end{equation}
where $(\cdot, \cdot)$ is $L^2$ inner product, and $a(u,v) = \e (u',
v') + b(u,v') + c(u,v)$.

For a positive integer $n\ge 2$, let $\T^n$ be  an arbitrary grid of
the form
\begin{equation} \label{grid1}
\T^n = \{x_i| 0 = x_0<x_1<\cdots<x_{n+1} = 1\},
\end{equation}
and let $h_i = x_i - x_{i-1}$. By $\phi_i (x)$, we denote the nodal
basis function at $x_i$ for $1\le i \le N$ by
\begin{equation} \label{base1}
\phi_i(x) =
\begin{cases}
 \disp{\frac {x-x_{i-1}}{h_i}} \textrm{ if } & x\in [x_{i-1}, x_i] \\
 \disp{\frac{x_{i+1} - x}{h_{i+1}}} \textrm{ if } & x \in [x_i,
 x_{i+1}]\\
 0 & \textrm{ otherwise. }
\end{cases}
\end{equation}
The finite element space is defined by $V^n = \{ v^n | v^n =
\sum_{i=1}^n v^n_i \phi_i(x) \}$. The finite element discretization
of (\ref{wksol1}) is to find $u^n \in V^n$ such that
\begin{equation} \label{fem1}
a(u^n, v^n) = (f,v^n), \forall v^n \in V^n \cap H_0^1.
\end{equation}
Existence and uniqueness of $u^n$ can be found in \cite{Brenner}
and  references therein. Now we denote
\begin{equation} \label{fem2}
u^n = \sum_{i=1}^n u^n_i \phi_i.
\end{equation}
Rewrite (\ref{fem1}) as
\begin{equation} \label{fem3}
\sum_{i=1}^n u^n_i a(\phi_i, \phi_j) = (f, \phi_j), j = 1, 2, \ldots
n.
\end{equation}
Let $A$ be an $n\times n$  matrix with
%$i$th row and $j$th column filled by
\begin{equation} \label{mat1}
a_{ij} = a(\phi_j,\phi_i) .
\end{equation}
Detailed calculation leads to further specific form of
\begin{equation} \label{mat2}
%\bdd
\barray
 a_{i,i} \ad = \e (\frac 1 {h_i} + \frac 1 {h_{i+1}}) + \frac c 3
(h_i + h_{i+1}) \\
 a_{i,i-1}\ad = - \frac \e {h_i} + \frac b 2 + \frac c 6 h_i\\
 a_{i, i+1} \ad= -\frac \e {h_{i+1}} - \frac b 2 + \frac c 6 h_i\\
 a_{i,j} \ad= 0, \quad \textrm{ if } |i-j| \ge 2 \earray
\end{equation}
Let $U^n = (u^n_1, \ldots, u^n_n)'$  and $F =
((f,\phi_1),\ldots, (f,\phi_n))'$ be  column vectors.
Then, (\ref{fem3}) is equivalent to  the linear system of
equations
\begin{equation}\label{fem4}
A U^n = F.
\end{equation}

Typically a FEM solution of a singularly perturbed BVP problem has boundary layer in a small interval  (associated with $\e$) of rapid variations of $u''$. Throughout this paper, unless it's explicitly mentioned, we assume solution $u$ of (\ref{bvp1}) has a boundary layer at $x=1$ and $x_n$ is located outside the boundary layer. This is reasonable assumption due to the very short interval of boundary layer depending on $0<\e\ll 1$. All
the results below can be obtained analogously for any boundary layer
located in
% between
$[0,1]$.

Let $\T^{n+m} = \T^n \cup \{s_1, \ldots, s_m\}$, where
$x_n<s_1<\cdots<s_m<x_{n+1}$. Denote the nodal basis functions on
$\T^{n+m}$ by $\{\phi_1, \ldots, \phi_{n-1}, \wdt \phi_n,
\phi_{s_1}, \ldots, \phi_{s_m}\}$, where $\wdt \phi_n$ and
$\phi_{s_i}$ are nodal basis for $x_n$ and $s_i$,
respectively.
% for each.
Note
that the first $n-1$ nodal basis functions of $\T^{n+m}$ are
exactly the same as those of $\T^n$.  Let $V^{n+m}$ be the function
space with basis $\{\phi_1, \ldots, \phi_{n-1}, \wdt \phi_n,
\phi_{s_1}, \ldots, \phi_{s_m}\}$. It is obvious that $V^n \subset
V^{n+m}$. Write $u^{n+m} $, the FEM solution of (\ref{bvp1}) in
$V^{n+m}$, as
\begin{equation} \label{wksol2}
u^{n+m} = \sum_{i=1}^{n-1} u^{n+m}_i \phi_i + u^{n+m}_n \wdt \phi_n
+ \sum_{i=1}^m u^{n+m}_{s_i}  \phi_{s_i}.
\end{equation}

In the next section, we fix $\T^n$, and start with observation
on the intersections of $u^n$ and $u^{n+m}$ for different
$\T^{n+m}$. For convenience, we use $Q_i \in u^n\cap u^{n+m}$ to denote the intersetion of $u^n$ and $u^{n+m}$ in the interval $(x_{i-1}, x_i)$, and by $x(Q_i)$ and $y(Q_i)$ we denote $x$- and $y$- coordinate of $Q_i$ respectively. The result shows that the intersections $\{Q_i: 2\le i \le n\}$ are independent of $m$ and distribution of $s_i$. Therefore, by adding
only one point $\{s_1\}$, we can compute $\{Q_i\in u^n\cap u^{n+1}\}$, and the accuracy of $Q_i$ has the same accuracy as $u^{n+\infty}$.

\section{Intersections of $u^n$ and $u^{n+m}$}
\begin{thm}\label{isc} \bdd {\rm
Fix $\T^n$. By adding one point $s_1\in (x_n, 1)$ arbitrarily,
we obtain new grid $\T^{n+1}$. Then the intersection $Q_i$ of $u^n$ and
$u^{n+1}$ in the interval $(x_{i-1}, x_{i})$  is independent of the choice of $s_1$ for any $i= 2, 3, \ldots, n$. That is, those coordinates of intersections do not
depend on the choice of $s_1\in (x_n, x_{n+1})$.

}
\end{thm}

\para{Proof.}
Analogous to (\ref{fem3}), we have a system of linear equations with
respect to $\{ u^{n+1}_i, i=1, \ldots, n; u^{n+1}_{s_1}\}$, given by
\begin{equation} \label{isc-1}
\sum_{i=1}^{n-1} u^{n+1}_i a(\phi_i, \phi_j) + u^{n+1}_n a(\wdt
\phi_n, \phi_j) + u^{n+1}_{s_1} a(\phi_{s_1}, \phi_j) = (f, \phi_j),
\quad j=1, 2,\ldots, n-1,
\end{equation}
\begin{equation} \label{isc-2}
\sum_{i=1}^{n-1} u^{n+1}_i a(\phi_i, \wdt \phi_n) +  u^{n+1}_n
a(\wdt \phi_n, \wdt \phi_n) + u^{n+1}_{s_1} a(\phi_{s_1}, \wdt
\phi_n) = (f, \wdt \phi_n),
\end{equation}
and
\begin{equation} \label{isc-3}
\sum_{i=1}^{n-1}  u^{n+1}_i a(\phi_i, \phi_{s_1}) + u^{n+1}_n a(\wdt
\phi_n,  \phi_{s_1}) +  u^{n+1}_{s_1} a(\phi_{s_1},  \phi_{s_1}) =
(f, \phi_{s_1}),
\end{equation}
Note that for $1\le j \le n-1$, $a(\wdt \phi_n, \phi_j) = a(\phi_n,
\phi_j)$ and $a(\phi_{s_1}, \phi_j) = 0$, and (\ref{isc-1}) leads to
\begin{equation}\label{isc-4}
\sum_{i=1}^{n}  u^{n+1}_i a(\phi_i, \phi_j) = (f, \phi_j), \quad
j=1, 2,\ldots, n-1.
\end{equation}
On the other hand, for $1\le i \le n-1$, $a(\phi_i, \wdt \phi_n) =
a(\phi_i, \phi_n)$, and (\ref{isc-2}) yields
\begin{equation} \label{isc-5}
\sum_{i=1}^{n-1} u^{n+1}_i a(\phi_i, \phi_n) + u^{n+1}_n a(\wdt
\phi_n, \wdt \phi_n) = (f, \wdt \phi_n) - u^{n+1}_{s_1} a(
\phi_{s_1}, \wdt \phi_n).
\end{equation}
For $1\le i \le n-1$, $a(\phi_i, \phi_{s_1}) = 0$, so it follows
from (\ref{isc-3})
\begin{equation} \label{isc-6}
u^{n+1}_n a(\wdt \phi_n, \phi_{s_1}) = (f, \phi_{s_1}) -
u^{n+1}_{s_1} a(\phi_{s_1}, \phi_{s_1}) .
\end{equation}
Let $p = (1-s)/h_{n+1}$. Observe $\phi_n = \wdt \phi_n +
p\phi_{s_1}$. Combining two equations above according to
(\ref{isc-5})$+ p * $(\ref{isc-6}), we have
\begin{equation} \label{isc-7}
\sum_{i=1}^{n-1} u^{n+1}_i a(\phi_i, \phi_n) + u^{n+1}_n a(\wdt
\phi_n,  \phi_n) = (f,  \phi_n) - u^{n+1}_{s_1} a(\phi_{s_1},
\phi_n).
\end{equation}
Hence,
\begin{equation} \label{isc-8}
\sum_{i=1}^{n} u^{n+1}_i a(\phi_i, \phi_n) = (f,  \phi_n) -
u^{n+1}_{s_1} a(\phi_{s_1}, \phi_n) + p u^{n+1}_n a(\phi_{s_1}
,\phi_n).
\end{equation}
Let $ U^{n+1} = (u^{n+1}_1, \ldots,u^{n+1}_n)'$ be a column vector
with length $n$.
By (\ref{isc-4}) and (\ref{isc-8}),
\begin{equation} \label{isc-9}
A U^{n+1} = \wdt F,
\end{equation}
where $\wdt F$ is a column vector with left-hand side of (\ref{isc-4})
and (\ref{isc-8}) as elements.
Subtracting  (\ref{isc-9}) from (\ref{fem4}),
\begin{equation} \label{isc-10}
%A (U - U^{n+1})
A (U^n - U^{n+1})
= F- \wdt F
\end{equation}
Notice that $F - \wdt F = C_{s,1} e_n$, where $e_n = (0,\ldots,0,1)'$ is
a vector with length $n$, and $C_{s,1} = u^{n+1}_{s_1} a(\phi_{s_1},
\phi_n) - p u^{n+1}_n a(\phi_{s_1}, \phi_n)$. Note that $C_{s,1}$ is a
scalar depending only on $s_1$, since $u_n^{n+1}$ term in $C_{s,1}$ is completely
determined by $s_1$. Therefore,
\begin{equation} \label{isc-11}
U^n -  U^{n+1} = C_{s,1} A^{-1} e_n
\end{equation}
The last equation tells us every $u^n_i - u^{n+1}_i$ increases or
decreases by the factor $C_{s,1}$ uniformly in $i$. Using elementary similar
triangle properties, we prove the result. \qed

\begin{rem} \label{isc-re} \bdd {\rm
If $u^n_i - u^{n+1}_i$ and $u^n_{i+1} - u^{n+1}_{i+1}$ have opposite sign, then $u^n$ and $u^{n+1}$ have intersection in $(x_i, x_{i+1})$. Notice that $A^{-1}e_n$ in (\ref{isc-11}) is FEM solution of green function of operator $A$. It is very common that FEM solution of green function intersects $x$-axis in each grid. Intuitively, this explains why  $u^n$ and $u^{n+1}$ intersect each other in every grid in most cases. Later we will present the criteria to be used for identifying the existence of intersections, see \lemref{int}.  Moreover, if there is no intersection in some interval $(x_i, x_{i+1})$ for a choice of $s_1$, then there will be no intersection for any choice of $s_1$.  }
\end{rem}

\begin{thm} \label{gsc} \bdd {\rm
Fix $\T^n$. Let $\T^{n+m}  = \T^n \cup \{s_1<s_2<\cdots <s_m\}$,
where $s_i \in (x_n, 1)$. Then the intersection $Q_i$ of $u^n$ and
$u^{n+m}$ in the interval $(x_{i-1}, x_{i})$ is independent of $m$ and distribution of $\{s_i\}$  for any $i=2, 3, \ldots, n$.
}
\end{thm}

\para{Proof.}
Let $ V^{n+m}$ be a function space with nodal basis functions
$\{\phi_1, \ldots, \phi_{n-1}, \wdt \phi_n,  \phi_{s_1}, \ldots,
\phi_{s_m} \}$ on $\T_{n+m}$. Analogous to (\ref{isc-4}), we have
\begin{equation}\label{gsc-1}
\sum_{i=1}^{n} u^{n+m}_i a(\phi_i, \phi_j) = (f, \phi_j), \quad j=1,
2,\ldots, n-1.
\end{equation}
Since $ V^{n+m} \supset V^n$, there exists a linear combination
$\phi_n = \wdt \phi_n  + \sum_{i=1}^m p_i \wdt \phi_{s_i}$ for some
$p_1, p_2, \ldots, p_m \in [0,1]$. Applying similar arguments
as that of
\thmref{isc}, we obtain
\begin{equation} \label{gsc-2}
\sum_{i=1}^{n}  u^{n+m}_i a(\phi_i, \phi_n) = (f,  \phi_n) +
\sum_{i=1}^m (p_i u^{n+m}_n - u^{n+m}_{s_i}) a( \phi_{s_i}, \phi_n).
\end{equation}
Define $C_{s,m} = \sum_{i=1}^m (u^{n+m}_{s_i} -p_i u^{n+m}_n) a(
\phi_{s_i}, \phi_n)$. Using exactly the same argument in
(\ref{isc-11}), we have
\begin{equation} \label{gsc-3}
U^n -  U^{n+m} = C_{s,m} A^{-1} e_n.
\end{equation}
Hence, the result follows. \qed

\begin{cor} \label{lsc} \bdd {\rm
Fix $\T^n$. Let $\T^{n+m}  = \T^n \cup \{s_1<s_2<\cdots <s_m\}$,
where $s_i \in (0, x_{1})$. Then the intersection $Q_i$ of $u^n$ and
$u^{n+m}$ in the interval $(x_{i-1}, x_i)$ is independent of $m$ and distribution of $\{s_i\}$ fixed for any $i= 2,\ldots, n$.  }
\end{cor}

\para{Proof.}
We rearrange the order of the index from $\{0, 1, 2, \ldots, n, n+1\}$ to
$\{n+1, n, \ldots, 1, 0\}$, and change the coordinate linearly from
$[0,1]$ into $[1,0]$. Using the same line of argument
as that of
\thmref{gsc}, the result holds.\qed

\begin{cor} \label{csc} \bdd {\rm
Fix $\T^n$. Let $\T^{n+m}  = \T^{n} \cup \{s_1<s_2<\cdots <s_m\}$,
where $s_i \in (x_{k-1}, x_{k})$ for some $2\le k \le n$. Then the
intersection of $u^n$ and $u^{n+m}$ in the interval $(x_{i-1}, x_{i})$
is independent of $m$ and distribution of $\{s_i\}$ for any
$i \in \{2, \ldots, n\} \setminus \{k\}$. }
\end{cor}
\para{Proof.}
This is straight forward result from \thmref{gsc} and \corref{lsc}.
\qed

\section{An $\e$-uniform Approximation $u^{n+1}$ in $[0,x_n]$}
In the previous section, by arbitrarily choosing a point $s_1 \in
(x_n,1)$, we can determine $Q_i\in u^n\cap u^{n+1}$ in each
interval, and the result shows $Q_i$ has the same accuracy
as that of
$u^{n+\infty}$. In this section,  by choosing
appropriate $\hat s_1 \in (x_n,1)$, we obtain $\hat u^{n+1}$, which
has $\e$-uniform accuracy in $[0,x_n]$. This will automatically
imply that $Q_i$ has $\e$-uniform accuracy, since $Q_i\in \hat
u^{n+1}$. For simplicity, we slightly abuse notation: Let $a_{n,\hat{s_i}} = a(\phi_{\hat{s_i}}, \phi_n)$ without confusing.

\begin{lem} \label{chs} \bdd {\rm
There exists $\hat s_1 \in (x_n, 1)$, such that, $a_{n,\hat s_1} =
0$ for $\hat \T^{n+1} = \{ x_0< x_1<\cdots<x_n<\hat s_1<x_{n+1}\}$.

}
\end{lem}

\para{Proof.}
By  (\ref{mat2}), to establish the desired result,
it is equivalent to prove that there exists
$0<h_{\hat s_1}<1-x_n$, satisfies
\begin{equation} \label{chs-1}
-\frac {\e} {\hat s_1} - \frac b 2 + \frac c 6 h_{\hat s_1} = 0,
\end{equation}
where $h_{\hat s_1} = \hat s_1 - x_n$. By eliminating the denominators in
the equation (\ref{chs-1}), we have
\begin{equation} \label{chs-2}
c h_{\hat s_1}^2 - 3b h_{\hat s_1} - 6\e = 0.
\end{equation}
If $c=0$, then $b<0$, and $h_{\hat s_1} = \disp{\frac {-2\e} b} >0$.
If $c\neq 0$, then the determinant of (\ref{chs-2}) is $9b^2 + 24 \e
c >0$. Write $h_{\hat s_1}$ using quadratic formula,
\begin{equation} \label{chs-3}
0< h_{\hat s_1} = \frac {3b + \sqrt{9b^2 + 24 \e c}}{2c} \le
\sqrt{\frac{6\e}{c}}.
\end{equation}
Thus, $h_{\hat s_1} = O(\e)$ if $c=0$, and $h_{\hat s_1} =
O(\sqrt{\e})$ if $c\neq 0$. \qed

\begin{rem} \label{remchs}\bdd{\rm
The essence of \lemref{chs} is to find such a $h_{\hat s_1}$ with $a_{n,\hat{s_1}} = 0$. If $b$ and $c$ are not constant, we can compute the formula for $h_{\hat s_1}$ involved with integrals. It is also possible to find it by discretizations.

}
\end{rem}

\begin{thm} \label{sth}\bdd{\rm
Given $\T^n$, take $\hat \T^{n+1}$ and $\hat s_1$ as in
\lemref{chs}. Use $\hat u^{n+1}$ to denote the FEM solution on
$\T^{n+1}$ of (\ref{bvp1}). Consider another BVP problem
\begin{equation} \label{sth-1}
-\e w'' - b w' + c w = f, \quad w(0) = 0, w(\hat s_1) = u(\hat s_1),
\end{equation}
where $u\cd$ is solution of (\ref{bvp1}). Use $w^n$ to denote the FEM
solution of (\ref{sth-1}) on $\hat \T^{n+1} \setminus \{1\}$, then
\begin{equation} \label{sth-2}
w(x) = u(x), \quad \forall x\in [0,\hat s_1],
\end{equation}
and
\begin{equation} \label{sth-3}
\hat u^{n+1} (x) = w^n(x), \quad \forall x\in [0,x_n].
\end{equation}

}
\end{thm}

\para{Proof.}
Note that $(\hat u^{n+1}_1, \hat u^{n+1}_2, \ldots, \hat u^{n+1}_n)$
is a
solution of the system of linear equations
\begin{equation} \label{sth-4}
\begin{cases} a_{i,i-1} \hat u^{n+1}_{i-1} + a_{i,i} \hat u^{n+1}_{i} + a_{i,i+1}
\hat u^{n+1}_{i+1} = (f, \phi_i) \quad i=1,2,\ldots,n-1 \\
a_{n,n-1} \hat u^{n+1}_{n-1} + a_{n,n} \hat u^{n+1}_n = (f, \wdt
\phi_n) -a_{n,\hat s_1} \hat u^{n+1}_{\hat s_1}.
\end{cases}
\end{equation}
Let $w^n = \sum_{i=1}^{n-1} w^n_i \phi_i + w^n_n \wdt \phi_n +
w^n_{\hat s_1} \phi^-_{\hat s_1}$, where $\phi^-_{\hat s_1} =
\phi_{\hat s_1}|_{[0,\hat s_1]}$. Then $(w_i^n$ for $i\in
\{1,2,\ldots, n, \hat s_1\})$ is a solution of the system of linear
equations
\begin{equation} \label{sth-5}
\begin{cases} a_{i,i-1} w^n_{i-1} + a_{i,i} w^n_{i} + a_{i,i+1}
w^n_{i+1} = (f, \phi_i) \quad i=1,2,\ldots,n-1 \\
a_{n,n-1} w^n_{n-1} + a_{n,n} w^n_n = (f, \wdt \phi_n) -a_{n,\hat
s_1}
w^n_{\hat s_1}\\
w_{\hat s_1}^n = u(\hat s_1) .
\end{cases}
\end{equation}
The solutions of
(\ref{sth-4}) and (\ref{sth-5}) are precisely the same, since
$a_{n,\hat s_1}=0$. \qed

\begin{rem}\label{srk} \bdd{\rm
From \thmref{sth}, we can separate the boundary layer by adding
point $\hat s_1\in (x_n,1)$. Therefore, it is equivalent to solve
 non-singularly perturbed
BVP problems by the FEM, and all general FEM error analysis
works well without effected by boundary layer. For example, if
$\{0<x_1<\cdots<x_n\}$ is uniform mesh in $[0,x_n]$, then
%maximum norm
$\|\hat u^{n+1}\|_{\infty,[0,x_n]}$ is bounded by
$\|u''\|_{\infty,[0,x_n+O(\sqrt\e)]} h^2$, and
$\|u''\|_{\infty,[0,x_n+O(\sqrt\e)]}$ is $\e$-uniformly bounded. On
the other hand, add $m$ points  in $(\hat s_1, 1)$, denoted by $\hat
\T^{n+1+m}$. Use $\hat \T^{n+1+\infty}$ to denote the grid which is
almost dense in $[\hat s_1,1]$. Use $\hat u^{n+1+\infty}$ to denote
the FEM solution of (\ref{bvp1}) on $\hat \T^{n+1+\infty}$. Then, $\hat
u^{n+1}$ is exactly the same with $\hat u^{n+1+\infty}$ on
$[0,x_n]$.

}
\end{rem}

\section{Auxiliary results}
%\subsection{Location of $Q_i$}
Recall $A$ is an $n\times n$ matrix with $a_{ij} = a(\phi_j,
\phi_i)$, and $e_n$ is $(0,\ldots,0,1)'$ of length $n$. Let $A_i$ be matrix replacing $i$th column of $A$ with $e_n$.
\begin{lem} \label{int} \bdd {\rm
Fix $\T^n$. Let $\T^{n+m}  = \T^n \cup \{s_1<s_2<\cdots <s_m\}$,
where $s_i \in (x_n, x_{n+1})$. Then $u^n$ and $u^{n+m}$ have their
intersection $Q_i$ in the interval $(x_{i-1}, x_{i})$ for some $2\le i
\le n$ if and only if
\begin{equation} \label{asp-2}
\det [A_{i} A_{i-1}] < 0,
\end{equation}
and the coordinates of $Q_i$ is given by
\begin{equation}\label{asp-1}
Q_i = \left(\frac {r_i} {r_i +1} x_{i} + \frac 1 {r_i +1} x_{i-1},
\frac {r_i} {r_i +1} u^n_{i} + \frac 1 {r_i +1} u^n_{i-1} \right) ,
\end{equation}
where $r_i = |\det A_{i-1} / \det A_{i}|$.
%Analogously,
%$u^n$ and $u^{n+m}$ have their intersection $Q_i$ in the interval
%$(x_i, x_{i+1})$ for some $k+1\le i \le n$ if and only if
%\begin{equation} \label{asp-3}
%\det [\wdt A_{i}^1 \wdt A_{i}^{2}] < 0,
%\end{equation}
%$Q_i$ is given by the same formula (\ref{asp-1}) with $r_i = |\det
%\wdt A_{i+1}^2 / \det \wdt A_{i+1}^{1}|$.
}
\end{lem}
\para{Proof.}
To obtain $Q_i$, we apply $\T^n$ and $\T^{n+m}$ to \thmref{gsc}.
Using Crammer's rule in (\ref{gsc-3}), we obtain
\begin{equation} \label{int-1}
u^n_i -  u^{n+m}_i = C_{s,m} \frac{\det A_i}{\det A}, \quad i = 1, 2,
\ldots, n ,
\end{equation}
Therefore
\begin{equation} \label{int-2}
\frac{u^n_i -  u^{n+m}_i}{u^n_{i-1} - u^{n+m}_{i-1}} =  \frac {\det
A_i}{\det A_{i-1}}, \quad i = 2, \ldots, n.
\end{equation}
A necessary and sufficient condition to have
an intersection is
$(u^n_i -  u^{n+m}_i)/(u^n_{i+1} - u^{n+m}_{i+1})<0$. This proves
(\ref{asp-2}). Using similar
%triangular property,
triangles, (\ref{asp-1})
follows.
%The counterpart for $i\ge k+1$ is symmetric proof.
%The converse can be proved in a similar way.
\qed

It is very common to have oscillation in finite element solution,
and we can use \lemref{int} to verify its behavior, see \remref{isc-re}. The following
theorem is a direct consequence of using Shishkin mesh.

\begin{thm}\label{skt} \bdd {\rm
Assume $\T^n$ is a uniform grid in $[0,1]$ satisfies condition
(\ref{asp-2}), and the boundary layer is at $x=1$. Then
\begin{equation} \label{skt-1}
\max_{1\le i \le n-1} | u^n(x(Q_i)) - u^n_I(x(Q_i))| < C n^{-2},
\end{equation}
where $C$ is independent of $\e$.

}
\end{thm}
\para{Proof.}
We put $m = O(n)$ grid in $(x_n,1)$, so that $\T^{n+m}$ forms
Bakhvalov grid or Shishkin grid. The uniform convergence of
$u^{n+m}$ on $\T^{n+m}$ is well known (see \cite{Kop1, Lin, Miller,
Zhang}) as
\begin{equation} \label{skt-2}
\|u^{n+m} - u_I^{n+m}\|_{\infty} \le C n^{-2}
\end{equation}
Also, we have $Q_i \in u^{n+m} \cap u^{n}$ by \corref{csc}. So
\begin{equation} \label{skt-3}
|u^{n+m}(x(Q_i)) - u_I^{n+m}(x(Q_i))| \le C n^{-2}.
\end{equation}
Note that $u_I^{n+m}|_{(0,x_n)} = u_I^{n}|_{(0,x_n)}$. Thus, the
theorem holds. \qed

\begin{rem} \bdd {\rm
From the result of \thmref{skt}, we have estimation of $O(h^2)$. In
non-uniform case, we can obtain
an error bound $O(h)$
% as error estimation
directly
from \cite{Chen}.

}
\end{rem}

\section{Numerical Results}
In this section, we present several examples.The first is a convection
diffusion equation, the second is a reaction diffusion equation, and
the
last one is a Green function.
\begin{exm} \bdd \label{cde} {\rm
Consider the convection-diffusion equation:
\begin{equation} \label{cde-1}
-\e u'' + u' = x, \quad u(0) = u(1) = 0.
\end{equation}
The exact solution is
\begin{equation} \label{cde-2}
u = x\left(\frac x 2 + \e \right) - \left(\frac 1 2 + \e\right)
\left(\frac {e^{(x-1)/\e} - e^{-1/\e}} {1 - e^{-1/\e}}\right).
\end{equation}
The solution $u\cd$ has a boundary layer at $x=1$, and is nearly
quadratic outside the boundary layer.

First, we use the linear finite element method on two different grid
$\T^{15}$ and $\T^{15+1}$ for $\e = 10^{-3}$, where $\T^{15}$ is a
uniform mesh on $[0,1]$ with $16$ intervals, and $\T^{15+1}$ is a
modified $\T^{15}$ with one point added at the center of the last
interval.  The intersections of finite element solution $u^{15}$ and
$u^{15+1}$ are almost on the interpolation of exact solution
$u_I^{15+1}$, as shown in Figure \ref{figcde}-a.

Second, we use the grid $\hat \T^{15+1}$ to compute for the same
$\e$, where $\hat \T^{15+1}$ is modified from $\T^{15}$ by adding
one specific point $\hat{s_1} \in (x_n,1)$ with $\hat{s_1} - x_n = 2\e$,
see \lemref{chs}. The finite element solution $\hat{u}^{15+1}$ is almost
overlapped with interpolation of interpolation of exact solution
$u_I^{15}$ in $[0,x_n]$, as seen from Figure \ref{figcde}-b. This verifies \thmref{sth}.

\begin{figure}[htbp]
\begin{center}
\mbox{ \subfigure[$u^{15}$ and $u^{15+1}$(dotted lines);
$u^{15+1}_I$(solid line)
]{\epsfig{figure=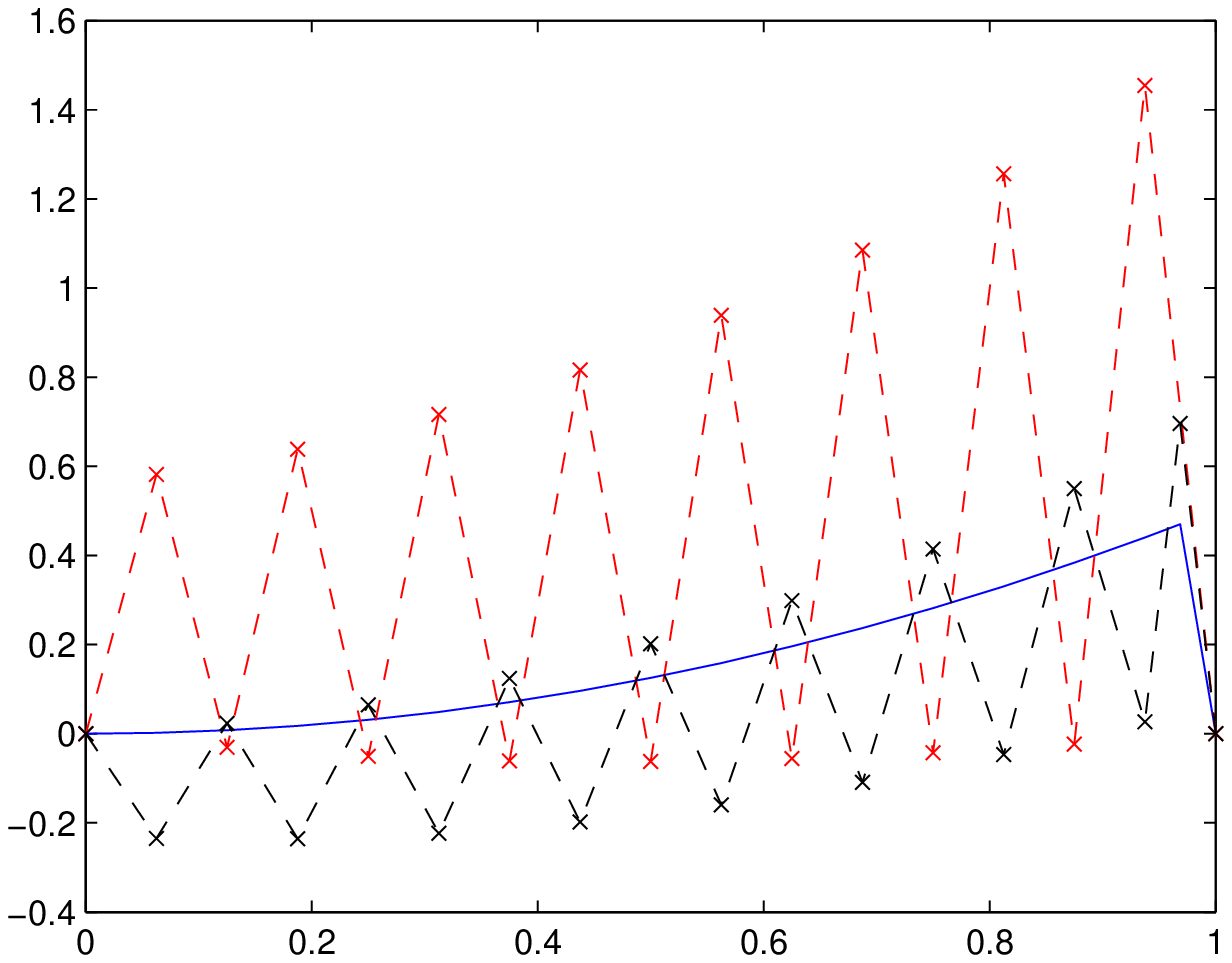,width=0.48\linewidth }}} \mbox{
\subfigure[$\hat u^{15+1}$(dotted line) and $u^{15}_I$ (solid line)
]{\epsfig{figure=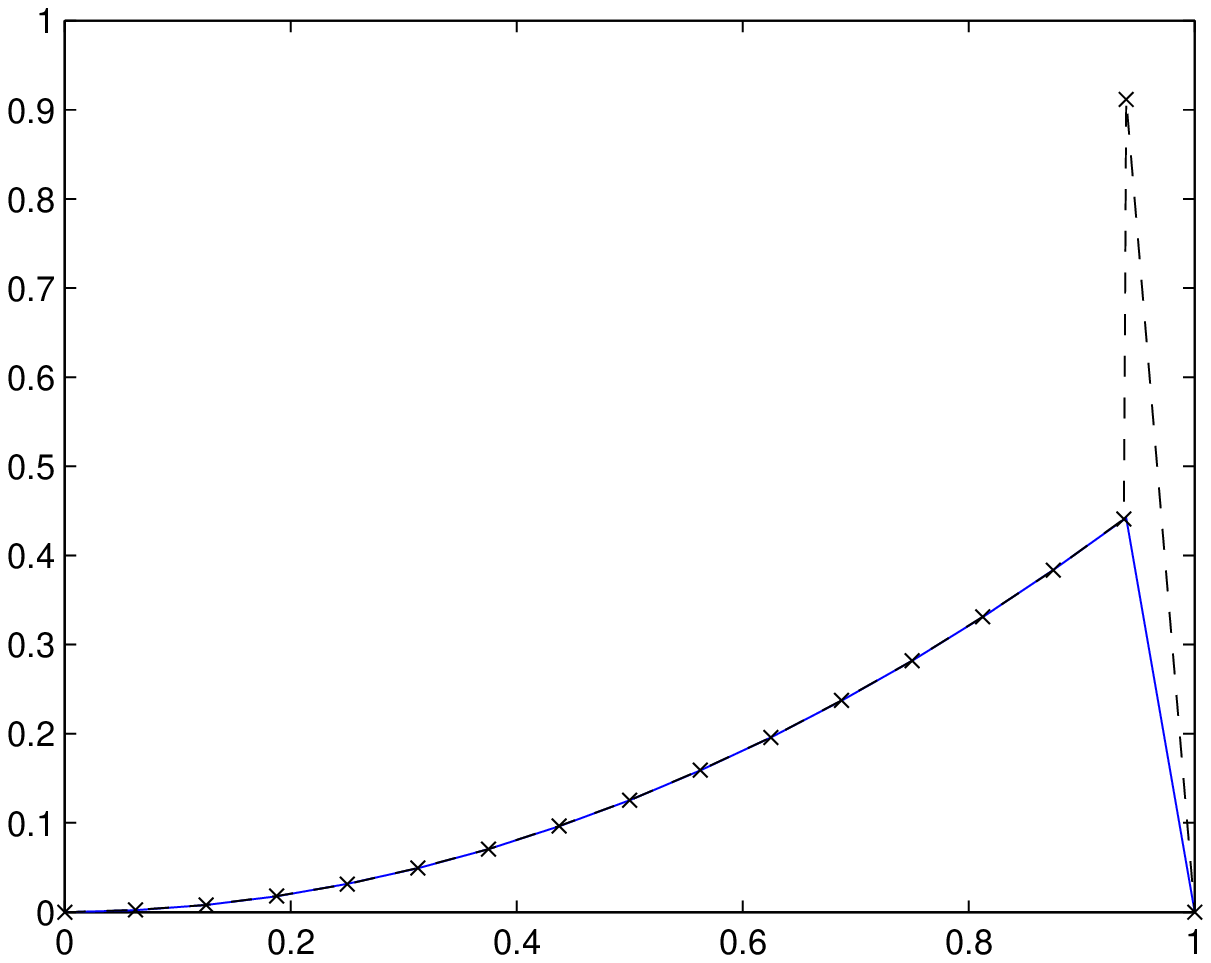,width=0.48\linewidth }}} \caption{
FEMs with $\e=10^{-3}$ for \exmref{cde} } \label{figcde}
\end{center}
\end{figure}

To compare with the well-known Shishkin mesh, we construct
$\T_s^{n+n}$, which divides both $[0,1-\theta]$ and $[1-\theta, 1]$ into $n$ equidistant subintervals, where $\theta = \{\frac 1 2, \frac{2\e\ln 2n}{b}\}$. $u_s^{n+n}$ is used to denote the FEM solution on $\T_s^{n+n}$.  Table \ref{tabcde1} shows the maximum norm of
$u_I^{n+1} - \hat u^{n+1}$ and $u_I^{n+1} - u_s^{n+n}$ in $[0,x_n]$.
Apparently, both $\hat u^{n+1}$ and $u_s^{n+n}$ has $\e$-uniform accuracy. However, $\hat u^{n+1}$ has better accuracy than $u_s^{n+n}$ by using less grids.
The reason is that $\hat u^{n+1}$ is
completely isolated  from the impact of errors from boundary layer; see
Table \ref{tabcde1}. This also verifies \thmref{sth}.

\begin{table}[htbp]
\begin{center}
\begin{tabular}{cccccc}
\hline
 & \multicolumn{2}{c}{$\e=10^{-5}$} & &\multicolumn{2}{c}{$\e=10^{-10}$}\\
\cline{2-3} \cline{5-6} n & $\|u_I^{n+1} - \hat u^{n+1}\|_{\infty,
[0,x_n]}$ & $\|u_I^{n+1} - u_s^{n+n}\|_{\infty, [0,x_n]}$
& & $\|u_I^{n+1} - \hat u^{n+1}\|_{\infty, [0,x_n]}$ & $\|u_I^{n+1} - u_s^{n+n}\|_{\infty, [0,x_n]}$\\
\hline 4 & 6.663e-003 &1.117e-002 & & 6.667e-003 & 1.117e-002\\
8 & 2.054e-003 &1.567e-003 & & 2.058e-003 & 1.569e-003\\
16 & 5.734e-004 &3.480e-004 & & 5.767e-004 & 3.500e-004\\
32 & 1.498e-004 & 8.384e-005& & 1.530e-004 & 8.569e-005\\
64 & 3.637e-005 & 1.948e-005& & 3.941e-005 & 2.115e-005\\
128 & 7.569e-006 & 3.928e-006& & 9.974e-006 & 5.221e-006\\
256 & 1.340e-006 & 1.340e-006& & 2.482e-006 & 1.292e-006\\
512 & 3.102e-007 & 6.738e-007& & 5.919e-007 & 3.208e-007\\
\hline
\end{tabular}
\end{center}
\caption{$\hat u^{n+1}$ and $u_s^{n+n}$ are FEM solutions on $\hat
\T^{n+1}$ and Shishkin mesh $\T_s^{n+n}$ for \exmref{cde}.
}\label{tabcde1}
\end{table}

Let $\e=10^{-10}$. Table \ref{tabcde2} shows the accuracy of $Q_i$,
the intersections of $u^8$ and $u^{8+1}$. Denote $x$- and $y$-
coordinates of $Q_i$ by $x(Q_i)$ and $y(Q_i)$, respectively.
Note that $Q_i$
has better accuracy than $\hat u^{8+1}$. The reason is yet to
%still open to
be discovered; see Table \ref{tabcde2}.

\begin{table}[htbp]
\begin{center}
\begin{tabular}{cccc}
    \hline
    $i$  &   $x(Q_i)$    &   $|y(Q_i) - u(x(Q_i))|$  & $|y(Q_i) - u^8_I(x(Q_i))|$\\
    \hline
       2   &     0.2499999996000000  &    7.499999579718697e-011  &4.999999719812465e-011\\
       3   &     0.2500000004000000  &     2.500008533523612e-011 & 8.326672684688674e-017\\
       4   &     0.4999999992000000  &    3.500000012035542e-010  &2.999999970665357e-010\\
       5   &     0.5000000008000000  &     5.000011515932101e-011 & 1.110223024625157e-016\\
       6   &     0.7499999988000000  &     6.625580639685325e-009 & 6.700580590379701e-009\\
       7   &     0.7500000012000000  &     7.500006171667906e-011 & 1.110223024625157e-016\\
    \hline
\end{tabular}
\end{center}
\caption{errors at $Q_i$ with $\e=10^{-10}$ on $\T^{8}$ and
$\T^{8+1}$ for  \exmref{cde}.}\label{tabcde2}
\end{table}

Plotted in Figure \ref{figcde1} are the convergence curves in the
maximum norm $\|u^{n+1}_I - \hat u^{n+1}\|_{\infty, [0,x_n]}$ for
$\e=10^{-5}$ and $\e = 10^{-10}$, respectively. They clearly indicate
the convergence rate is proportional to $n^{-2}$. It verifies Remark
\ref{srk}; see Figure \ref{figcde1}.

\begin{figure}[htbp]
\centering
\includegraphics[width=4in, height=4in]{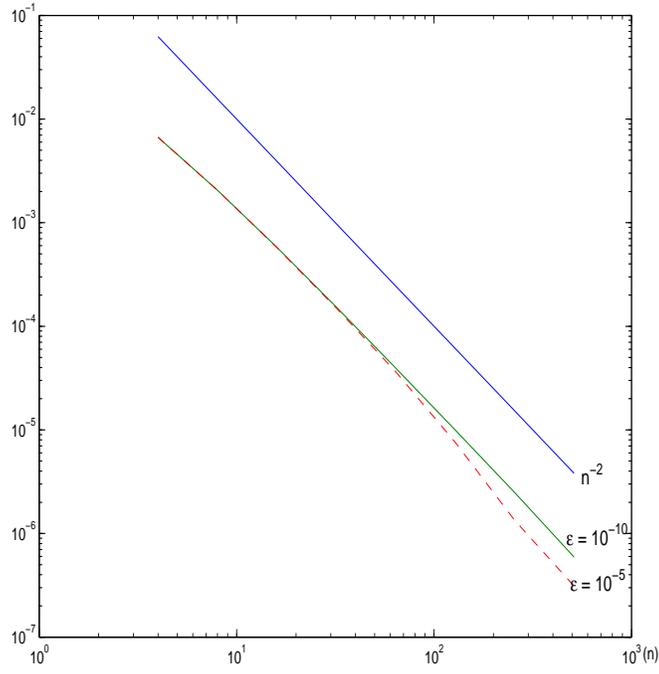}
\caption{errors $\|u - \hat u^{n+1}\|_{\infty, [0,x_n]}$ with
various $n$ for \exmref{cde}.}\label{figcde1}
\end{figure}

 }
\end{exm}

\begin{exm} \label{rde} \bdd{\rm
We examine the problem of
a reaction diffusion equation as another example of
(\ref{bvp1}).
\begin{equation} \label{bvp2}
-\e u''(x)+u(x)=x, \quad u(0) = u(1) = 0.
\end{equation}
The exact solution is
\begin{equation}\label{sol1}
u(x) = x-\frac{e^{(x-1)/\sqrt{\e}}-e^{-(x+1)/\sqrt{\e}}}{1-e^{-2/\sqrt{\e}}}
\end{equation}

The exact solution $u\cd$ has boundary layer at $x=1$, and is nearly
linear outside the boundary layer. Also, reaction diffusion equation has
relatively stable matrix $A$ compared with convection diffusion
equation. Due to these reasons, the FEM solutions of
(\ref{bvp2}) is better than the FEM solutions of (\ref{cde-1}).

For $\e = 10^{-10}$, we compute
the FEM solution $u^4$ and $u^{4+1}$
on the grid $\T^{4}$ and $\T^{4+1}$, where $\T^4$ is uniform mesh on
$[0,1]$ and $\T^{4+1}$ is modified by adding one point at the center
of last interval; Figure \ref{figrde1}-a.

By adding one point $\hat{s_1}\in (x_n, 1)$ with $\hat{s_1} - x_n = \sqrt[]{6\e}$ as in \lemref{chs}, we use new grid $\hat{\T}^{n+1}$, and denote its FEM solution as $\hat{u}^{4+1}$. As shown in Figure \ref{figrde1}-b, $\hat u^{4+1}$ is almost overlapped with $u_I^4$, the interpolation of exact solution; see Figure \ref{figrde1}-b.

\begin{figure}[htbp]
\begin{center}
\mbox{ \subfigure[$u^4$ and $u^{4+1}$(dotted lines);
$u_I^{4+1}$(solid line)]
{\epsfig{figure=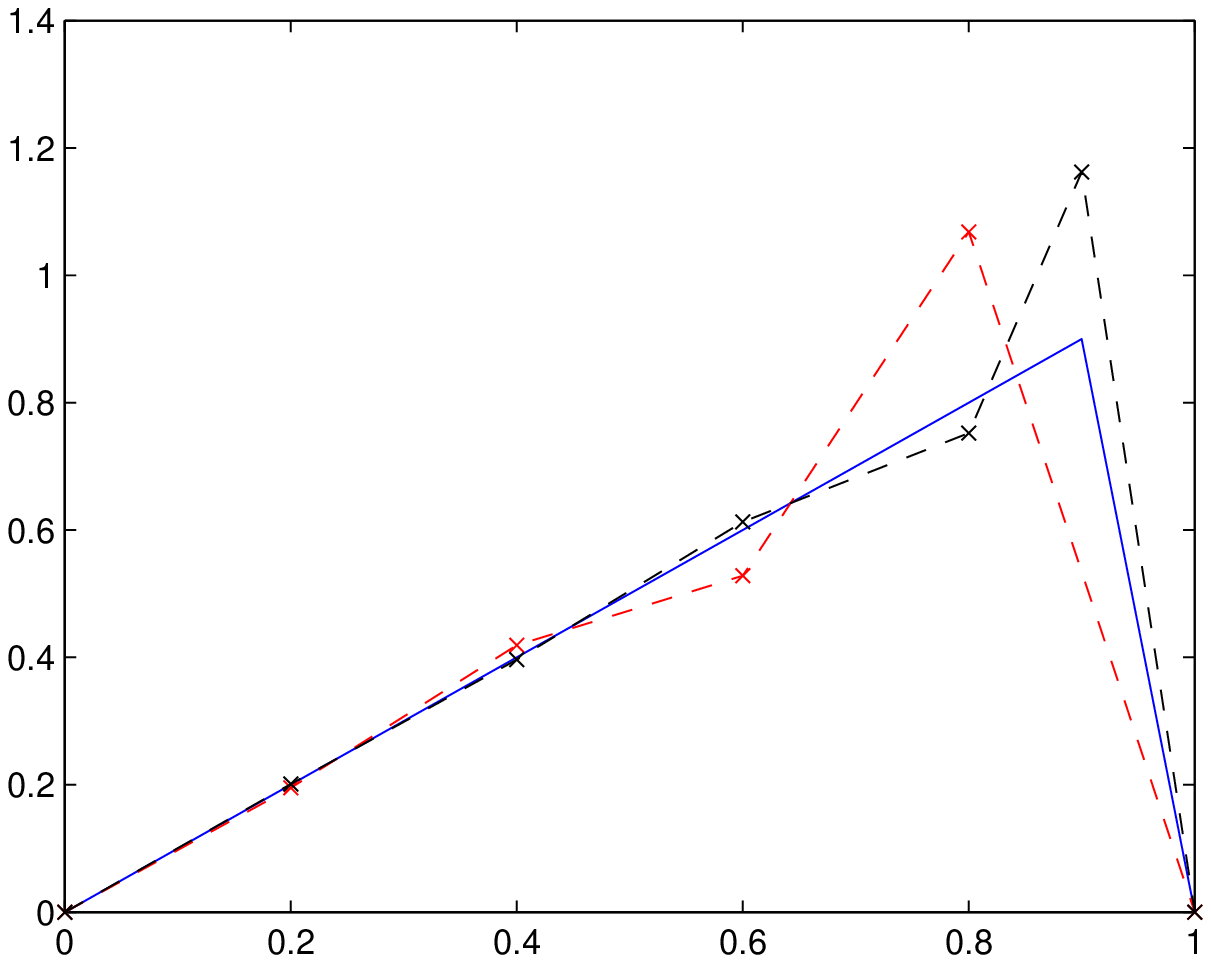,width=0.48\linewidth }}} \mbox{
\subfigure[$\hat u^{4+1}$(dotted line) and $u^4_I$(solid line)
]{\epsfig{figure=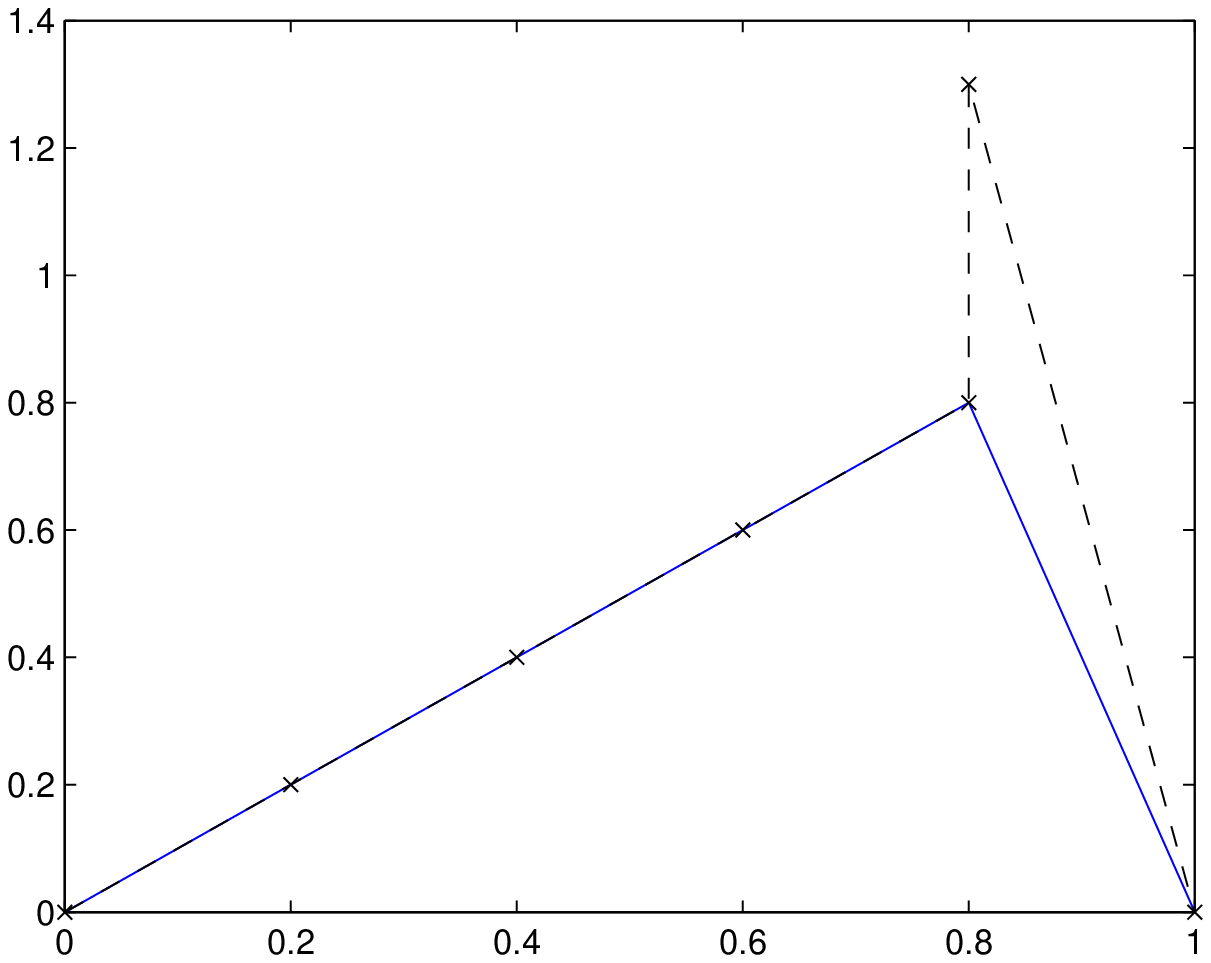,width=0.48\linewidth }}} \caption{FEMs
with $\e=10^{-10}$  for \exmref{rde}. } \label{figrde1}
\end{center}
\end{figure}

Let $\theta = \min\{\frac 1 2, \frac{\sqrt[]{\e} \ln 2n}{\sqrt{c}}\}$. We construct shishkin mesh $\T_s^{n+n}$ by dividing $[0,1-\theta]$ and $[1-\theta, 1]$ into $n$ equidistant subintervals.   Table \ref{tabrde1} present the errors of $\hat u^{n+1}$. Compared
with $u_s^{n+n}$, the FEM solutions using Shishkin mesh $\T_s^{n+n}$,
the errors are smaller and $\e$-uniform. We omit the convergence
curve and error table of $Q_i$, since all those errors are within
computer errors (around $10^{-14}$).

\begin{table}[htbp]
\begin{center}
\begin{tabular}{cccccc}
\hline
 & \multicolumn{2}{c}{$\e=10^{-5}$} & &\multicolumn{2}{c}{$\e=10^{-10}$}\\
\cline{2-3} \cline{5-6} n & $\|u_I^{n+1} - \hat u^{n+1}\|_{\infty,
[0,x_n]}$ & $\|u_I^{n+1} - u_s^{n+n}\|_{\infty, [0,x_n]}$
& & $\|u_I^{n+1} - \hat u^{n+1}\|_{\infty, [0,x_n]}$ & $\|u_I^{n+1} - u_s^{n+n}\|_{\infty, [0,x_n]}$\\
\hline
4 & 1.665e-016 &1.517e-004 &  &1.110e-016  & 4.980e-007\\
8 & 1.110e-016 & 5.415e-005 &  &2.220e-016  & 1.868e-007\\
16 & 2.220e-016 &2.161e-005 &  &3.331e-016  & 8.451e-008\\
32 & 2.220e-016  & 7.391e-006& & 3.331e-016 & 4.054e-008\\
64 & 3.331e-016  & 1.300e-006& & 5.551e-016 & 1.984e-008\\
128 & 4.441e-016 & 1.159e-009&  & 5.551e-016 & 5.551e-016\\
256 & 2.459e-013  & 2.948e-007& & 6.661e-016 & 9.795e-009\\
512 & 5.440e-015 & 2.865e-007&  & 7.772e-016 & 4.847e-009\\
\hline
\end{tabular}
\end{center}
\caption{$\hat u^{n+1}$ and $u_s^{n+n}$ are
the FEM solutions on $\hat
\T^{n+1}$ and Shishkin mesh $\T_s^{n+n}$ for \exmref{rde}.
}\label{tabrde1}
\end{table}

}
\end{exm}

\begin{exm} \label{gre} \bdd{\rm
This example presents a demonstration of \corref{csc}.
Using the FEM, we aim to find the Green function (as a solution of)
\begin{equation} \label{green}
-\e^2 u'' + u = \delta_\al, \quad u(0) = 0, u(1) = 0,
\end{equation}
where $\delta_\al$ is delta function with peak at $\al\in (0,1)$.
Denote a function as
\begin{equation} \label{gr1}
g(x) = \frac {e^{x/\e} - e^{-x/\e}} {e^{\al/\e} - e^{-\al/\e}}.
\end{equation}
The exact solution of (\ref{green}) is
\begin{equation} \label{grsol}
u = \begin{cases} K_\al g(x) & 0\le x \le \al \\
K_\al g\left(\frac{\al (1-x)}{1-\al}\right) & \al < x \le 1
\end{cases}
\end{equation}
where $K_\al\simeq \al/\e$ is a constant depend on $\al$ and $\e$.

$\T^{6+1}$ used in Figure \ref{figgr}-a is modified from $\T^6$ by
adding one point $s_1$ at the center of the last interval, while
$\T^{6+1}$ used in Figure \ref{figgr}-b is modified by adding $s_1$
at the center of $(x_3,x_4)$.

\begin{figure}[htbp]
\begin{center}
\mbox{ \subfigure[$u^6$ and $u^{6+1}$ (dotted lines) for $\frac 6
7<\al<1$. ]{\epsfig{figure=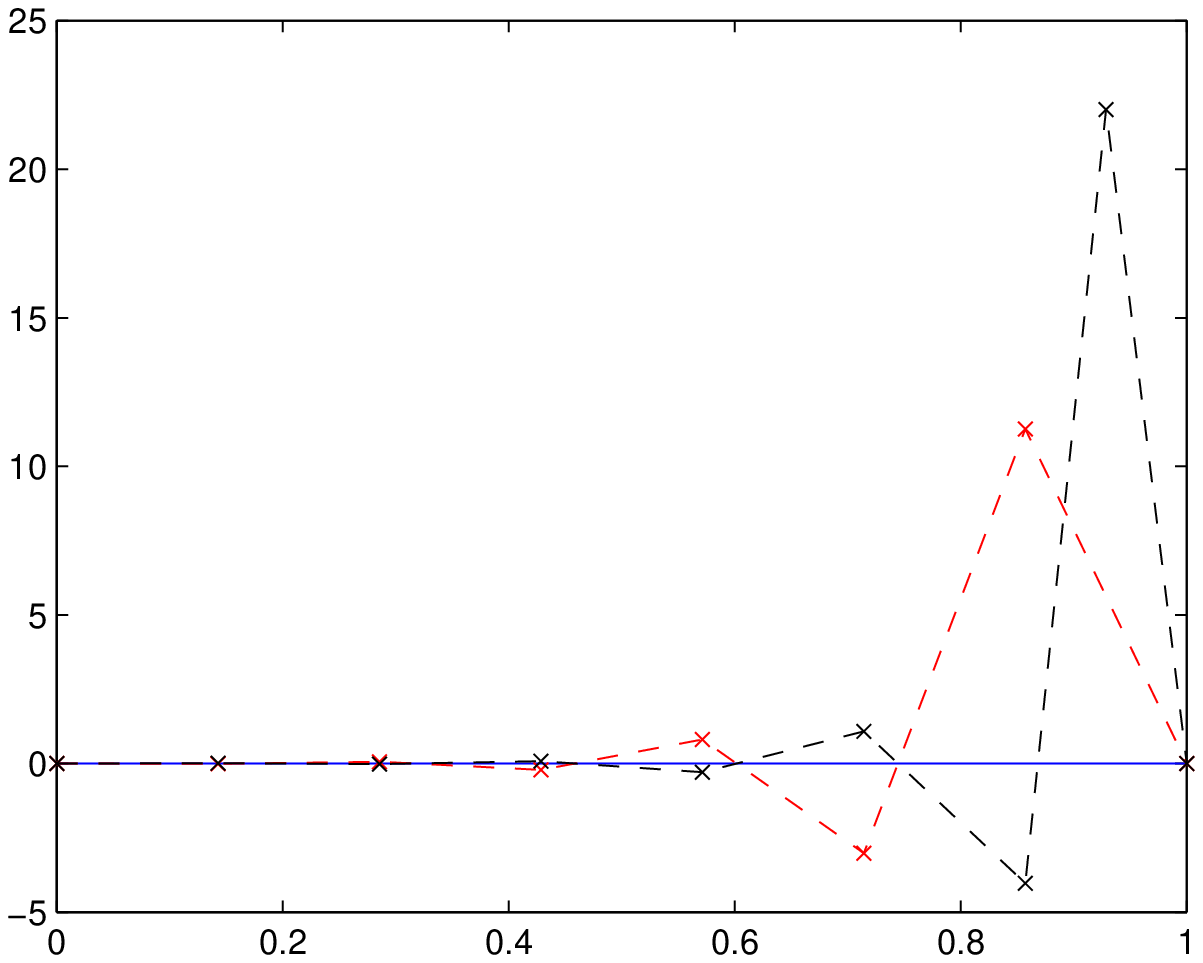,width=0.48\linewidth }}}
\mbox{ \subfigure[$u^6$ and $u^{6+1}$ (dotted lines) for $\frac 3
7<\al<\frac 4 7$. ]{\epsfig{figure=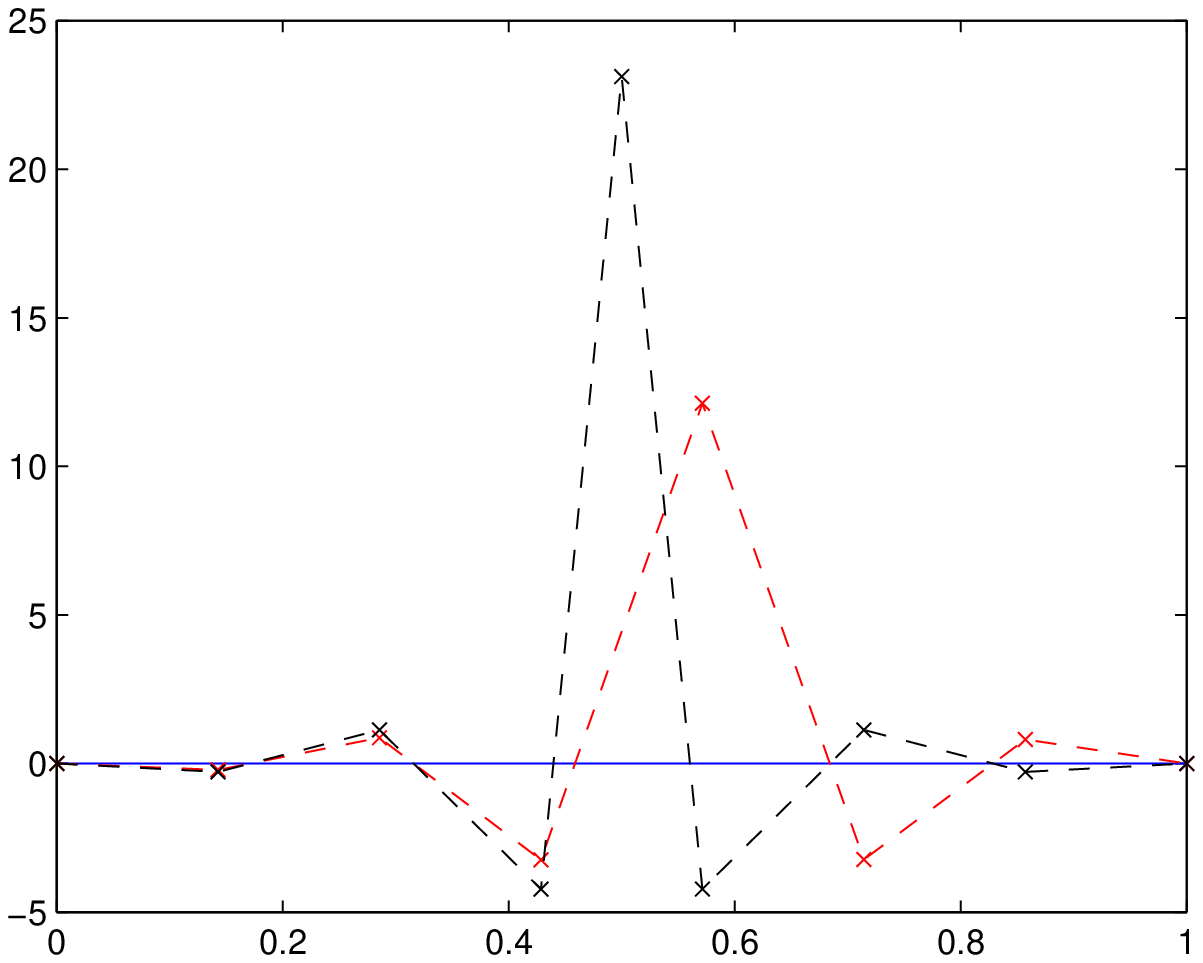,width=0.48\linewidth
}}} \caption{FEMs with $\e=10^{-5}$ for \exmref{gre}. } \label{figgre1}
\label{figgr}
\end{center}
\end{figure}

\begin{table}[htbp]
\begin{center}
\begin{tabular}{cccccc}
\hline
 & \multicolumn{2}{c}{$\frac 6 7 < \al < 1$} & &\multicolumn{2}{c}{$\frac 3 7 < \al < \frac 4 7$}\\
\cline{2-3} \cline{5-6} $i$ & $ x(Q_i) $ & $|y(Q_i) - u(x(Q_i))| $
& & $ x(Q_i) $ & $|y(Q_i) - u(x(Q_i))| $\\
\hline
2 & .1714   &    1.7347e-018 &  &.1714   &    2.7756e-017\\
3 &     .3158  &    6.9389e-018 &  &.3158   &    1.1102e-016\\
4 & .4588  &    2.7756e-017&  &--  & --\\
5 & .6016   &    0         & & .6842   &    1.7764e-015\\
6 & .7445     &  0         & & .8286  &    4.4409e-016\\
\hline
\end{tabular}
\end{center}
\caption{errors at $Q_i$ for $\e= 10^{-5}$ for \exmref{gre}.
}\label{tabgre1}
\end{table}

}
\end{exm}

\section{Further Remarks}
This paper is devoted to finite element methods for singularly
perturbed boundary value problems. An interesting behavior is
discovered: One can add arbitrary many points in one of the grids,
while the corresponding FEM solutions always have the common
intersections $\{Q_i\}$ in all other intervals. Moreover, a
practical and efficient $\e$-uniform mesh is developed. The FEM
solution under this mesh can be viewed as a non-singularly perturbed
BVP perturbation
problem, and all general FEM error analysis can be applied.

In both \exmref{rde} and \exmref{gre}, the errors are within
computer error. However, the errors of \exmref{cde} is visible
errors relative to computer error. The main reason is the exact
solution of \exmref{cde} is almost quadratic, while our
approximation is based on linear finite element space. To increase
accuracy, one can generalize the results to the higher order finite
element space. If the exact solution has several boundary layers, it
 can also be generalized to isolate each boundary layer.

Although the exact solution of \exmref{cde} is nearly quadratic, the
accuracy of intersections $\{Q_i\}$ is almost within computer error.
We know $\hat u^{n+1}$ has the accuracy of $\hat u^{n+1+\infty}$,
while $\{Q_i\}$ has the accuracy of $u^{n+\infty}$. The only
difference of two is the interval $(x_n, \hat s_1)$ of width
$O(\sqrt{\e})$ or $O(\sqrt{\e})$. In fact, this causes the error
difference from \exmref{cde}. It might be interesting to discover
the reason behind. It leads to the error analysis of
non-quasiuniform meshes.

It is very challenged to generalize the idea to isolate boundary
layer in higher dimensional cases. On the other hand, \lemref{int}
provided a necessary and sufficient condition to verify the behavior
of oscillation of specific FEM solution. However, it is not handy
enough to explain why the oscillation behavior is common to FEM
solutions. In general,
the problem of determining
%we still do not know
in what cases the FEM
solutions will or will not oscillate remains open.

\end{document}